\documentclass[letterpaper, 10 pt, conference]{ieeeconf}
 \IEEEoverridecommandlockouts
\usepackage{graphicx}
\usepackage{amsmath}
\usepackage{amssymb}

\newcommand{\figwidth}{72mm}
\newtheorem{theorem}{Theorem}

\newcommand{\ie}{\textit{i}.\textit{e}.}

\title{\LARGE \bf
	Sequential detection of Replay attacks
}

\author{Arunava Naha$^{1}$, André Teixeira$^{1}$, Anders Ahlén$^{1}$ and Subhrakanti Dey$^{2}$
	\thanks{*This work is supported by The Swedish Research Council (VR) under grants 2017-04053 and 2018-04396, and by the Swedish Foundation for Strategic Research.}
	\thanks{$^{1}$Arunava Naha, André Teixeira, and Anders Ahlén are with the Department of Electrical Engineering, Uppsala University, 751 03 Uppsala, Sweden
		{\tt\small arunava.naha@angstrom.uu.se, andre.teixeira@angstrom.uu.se, and Anders.Ahlen@angstrom.uu.se}}%
	\thanks{$^{2}$Subhrakanti Dey is with the Department of Electronic Engineering, Hamilton Institute, National University of Ireland, Maynooth, Ireland. He is also with the Department of Electrical Engineering, Uppsala University, 751 03 Uppsala, Sweden
		{\tt\small Subhra.Dey@signal.uu.se}}%
}

\begin{document}

\maketitle
\thispagestyle{empty}
\pagestyle{empty}

\begin{abstract}
One of the most studied forms of attacks on the cyber-physical systems is the replay attack. The statistical similarities of the replay signal and the true observations make the replay attack difficult to detect. In this paper, we have addressed the problem of replay attack detection by adding watermarking to the control inputs and then performed resilient detection using cumulative sum (CUSUM) test on the joint statistics of the innovation signal and the watermarking signal. We derive the expression of the Kullback-Liebler divergence (KLD) between the two joint distributions before and after the replay attack, which is asymptotically inversely proportional to the detection delay. We perform structural analysis of the derived KLD expression and suggest a technique to improve the KLD for the systems with relative degree greater than one. A scheme to find the optimal watermarking signal variance for a fixed increase in the control cost to maximize the KLD under the CUSUM test is presented. We provide various numerical simulation results to support our theory. The proposed method is also compared with a state-of-the-art method. 
\end{abstract}

\section{Introduction}
\label{sec:intro}
Nowadays, large cyber-physical systems (CPS) are getting deployed for intelligent transportation systems, manufacturing industries, smart grids, etc. \cite{Alguliyev2018}. Along with their immense advantages, there are also growing concerns about the safety and security of such systems. Attacks on the CPS can be a serious threat to the sensitive user data security, availability and reliability of critical resources, user's physical safety, and monetary loss \cite{Alguliyev2018}. Various techniques, such as data encryption,  authentication, firewall, cryptography, digital watermarking etc. are normally deployed to protect the cyber-layer of the CPS. Such protection schemes may not be adequate to protect the CPS from attacks on the physical layers as realised from different past incidents, such as the famous Stuxnet attack \cite{Langner2011}. In the Stuxnet attack, the malware issued harmful control inputs to increase the pressure of the centrifuges in a uranium enrichment plant in Iran \cite{Mo2015}. It also replaced the true measurements with previously recorded observations to remain stealthy. An attacker can launch a replay attack without detailed knowledge about the system parameters and control logic. The attacker can hijack a sensor node and record the observation for some time and then replay it back by replacing the true measurements at some later point of time. The attacker can alter the system in some harmful way and may remain stealthy during the replay attack.  

A widely applied technique for the replay attack detection is to add watermarking signal to the control inputs and then perform various statistical tests using the observations or the innovation signal from the Kalman estimator \cite{Mo2015, Satchidanandan2017, Ding2018a}. In one approach, $\chi^2$ statistics generated using the innovation signal is compared with some threshold for attack detection \cite{Zhao2020, Mo2015, Hosseinzadeh2019}. In another approach, test statistics are built using the observation data and performed a threshold check \cite{Satchidanandan2017} or the Neyman-Pearson (NP) test \cite{Zhai2020}. Addition of watermarking increases the probability of detection but at the same time, it increases the control cost \cite{Mo2015}. In \cite{Mo2015}, an optimal watermarking signal is designed which maximizes the attack detectability for a fixed increase in the control cost. In a different approach, the watermarking signal is also added or multiplied with the observations before the transmission. At the receiver, the authenticity of the observations are first checked, and then the watermarking signal is filtered out before feeding the observations to the estimator or controller \cite{Ferrari2017, Trapiello2019, Sanchez2019, Ye2019}. The watermarking signals for the observations could be of different types, such as sinusoidal \cite{Ferrari2017}, multiplicative to the observations \cite{Trapiello2019}, time-varying sinusoidal \cite{Sanchez2019}, random noise \cite{Ye2019}, etc. Since the added watermarking signal is removed before the observations are fed to the controller, such methods do not increase the control cost. However, if the attacker can access the signal before the addition of the watermarking, then these methods may fail. In \cite{Fang2020a}, the authors design a periodic watermarking scheme for the replay attack detection, which reduces the cost of adding the watermarking to the control inputs all the time before the attack. Even though most of the methods found in the literature studied the problem of replay attack detection for linear time-invariant (LTI) systems, a detection scheme is reported in \cite{Porter2020} for time-varying systems by adding time-varying dynamic watermarking.  There are few other methods found in the literature which do not use the watermarking for the replay attack detection. In \cite{Rath2020}, timestamps are added to the data, and in \cite{Hoehn2016}, a nonlinear element is inserted in the control loop for the replay attack detection. A set membership-based approach is followed in \cite{Liu2020}.  In addition to the research on replay attack detection, researchers have also studied the closed-loop stability of nonlinear systems under attack \cite{Huang2020}, the conditions on the watermarking to guarantee detection of the replay attack on dc microgrids \cite{Gallo2018}, and the state estimation problem when the system is under attack \cite{Franze2019, Chen2018b}.

Detection of an attack as early as possible is of immense importance for the CPS to reduce the magnitude of the damage. Most of the detection mechanisms reported in the literature do not address the issue of resilient detection of attacks explicitly. Moreover, some of the reported methods do batch processing which makes the detection delay dependent on the choice of the window size. In addition to that, since the processes are expected to run for a very long time before the attack takes place, the average run length (ARL) between the two false alarms is a better metric to use compared to the false alarm rate (FAR) \cite{Giraldo2019}. Therefore, we have applied a cumulative sum (CUSUM) test \cite{Tartakovsky2014, Girardin2018} using the joint distributions of the innovation signal and the watermarking signal before and after the replay attack. The added watermarking signal is independent and identically distributed (iid). In our previous paper \cite{watermarking_tac}, we have extensively studied two CUSUM tests, optimal CUSUM and sub-optimal CUSUM for the quickest detection of data deception attacks. In the data deception attack, the attacker generates fake observations using a linear stochastic process and replaces the true observations by them. In our current paper, we have modified the sub-optimal CUSUM test for the replay attack detection and studied its performance in terms of the average detection delay (ADD) and the increase in the control cost for a fixed upper bound on the ARL. We have derived the expression for the Kullback–Leibler divergence (KLD) measure between the joint distributions before and after the attack. KLD is asymptotically inversely proportional to the supremum of ADD (SADD) \cite{Tartakovsky2014, Girardin2018}. We have studied the effect of relative degree on the KLD and exhibited a way to improve the KLD for systems with relative degree greater than one. A technique to optimize the watermarking signal variance, which maximizes the KLD for a fixed upper bound on the increase in the control cost is also proposed.

The paper is organized as follows. The system model and the attack model are described in Section~\ref{sec:system_and_attack_model}. Section~\ref{sec:cusum} provides the replay attack detection scheme, the KLD expression for the replay attack, and the technique to improve the KLD for systems with a relative-degree greater than one. A technique for optimizing the watermarking signal variance is also discussed in Section~\ref{sec:cusum}. Section~\ref{sec:numerical_results} provides the numerical results and Section~\ref{sec:conclusion} concludes the paper. 

\section{System and Attack Model} \label{sec:system_and_attack_model}
This section discusses the system models during normal operations and under attack with the replay attack model considered in this paper. 
\subsection{System Model during Normal Operations}
\label{subsec:system_model_normal}
\begin{figure}[h!]
	\centering
	\includegraphics[width=45mm]{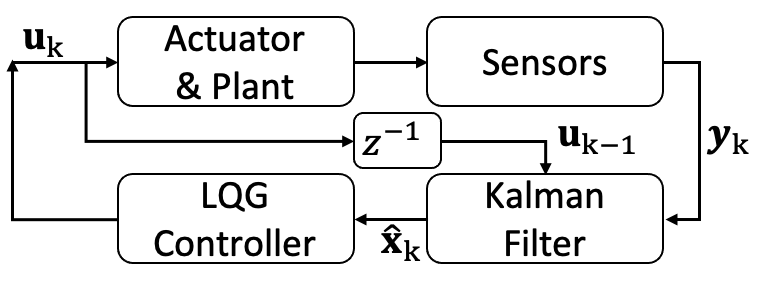}
	\caption{Schematic diagram of the system during normal operation.}
	\label{fig:sche_system_normal}
\end{figure}

Figure~\ref{fig:sche_system_normal} shows the schematic diagram of the network control system (NCS) under normal operation employed for this paper. The system is modelled as, 
\begin{align}
	&{\bf{x}}_{k+1}={\bf A}{\bf{x}}_{k}+{\bf B}{\bf{u}}_{k}+{\bf{w}}_{k},  \label{eqn:state_eqn} \ \\
	&{\bf{y}}_{k}={\bf C}{\bf{x}}_{k}+{\bf{v}}_{k}.
	\label{eqn:obj_eqn}
\end{align}
Here ${\bf{x}}_{k}\in {\rm I\!R}^{n}$, ${\bf{u}}_{k}\in {\rm I\!R}^{p}$, and ${\bf{y}}_{k} \in {\rm I\!R}^{m}$ are the state, input vector, and output vector  at the $k$-th time instant respectively.  ${\bf{w}}_{k} \in {\rm I\!R}^{n} \sim \mathcal{N}(0,{\bf Q})$ and ${\bf{v}}_{k} \in {\rm I\!R}^{m} \sim \mathcal{N}(0,{\bf R})$ are the iid process noise and observation noise respectively. ${\bf{A}}\in {\rm I\!R}^{n\times n}$, ${\bf{B}}\in {\rm I\!R}^{n\times p}$, ${\bf{C}}\in {\rm I\!R}^{m\times n}$, ${\bf{Q}}\in {\rm I\!R}^{n\times n}$ and ${\bf{R}}\in {\rm I\!R}^{m\times m}$. The noise vectors ${\bf{v}}_{k}$ and ${\bf{w}}_{k} $ are mutually independent, and both are independent of the initial state vector, ${\bf{x}}_{k_0}$. We assume the system is stabilizable and detectable. We also assume that the system has been operational for a long time, thus the system is currently at steady state. The states are estimated using a Kalman filer as follows,
\begin{align}
	{\hat{\bf{x}}}_{k|k-1} &={\bf A}{\hat{\bf{x}}}_{k-1|k-1}+{\bf B}{\bf{u}}_{k-1} \label{eqn:est_state_update}\ \\
	{\hat{\bf{x}}}_{k|k} &={\hat{\bf{x}}}_{k|k-1}+{\bf K}\gamma_k
	\label{eqn:kf_state_eqn}
\end{align}
where ${\hat{\bf{x}}}_{k|k-1} =E[{\bf{x}}_{k}|\Psi_{k-1}]$ and ${\hat{\bf{x}}}_{k|k} =E[{\bf{x}}_{k}|\Psi_{k}]$ are the predicted and filtered state estimates respectively. $E[\cdot]$ denotes the expected value and $\Psi_{k}$ is the set of all measurements up to time $k$. The innovation $\gamma_k$ and steady state Kalman filter gain $\bf K$ are given by 
\begin{align}
	\gamma_k&={\bf y}_k-{\bf C}{\hat{\bf{x}}}_{k|k-1} \label{eqn:gamma_k} \ \\
	{\bf K}&={\bf P}{\bf C}^T\left({\bf C}{\bf P}{\bf C}^T+{\bf R}\right)^{-1} \label{eqn:kalman_gain},
\end{align}
where ${\bf P}=E\left[({\bf x}_k-{\hat{ \bf x}}_{k|k-1}) ({\bf x}_k-{\hat{ \bf x}}_{k|k-1})^T   \right]$ is the steady state error covariance. ${\bf P} $ is the solution of the following algebraic Riccati equation
\begin{align}
	{\bf P}={\bf A}{\bf P}{\bf A}^T+{\bf Q}-{\bf A}{\bf P}{\bf C}^T\left( {\bf C}{\bf P}{\bf C}^T +{\bf R} \right)^{-1}{\bf C}{\bf P}{\bf A}^T.
	\label{eqn:P}
\end{align}
The control input ${\bf u}_k$ is generated by minimizing the infinite horizon LQG cost as given in \cite{Mo2015}. The optimal ${\bf u}_k^*$ takes the following form,
\begin{align}
	{\bf u}^*_k&={\bf L}{\hat {\bf x}}_{k|k}  \label{eqn:opt_u}, \ \\
	{\bf L}&=-\left( {\bf B}^T{\bf S}{\bf B}+{\bf U}\right)^{-1}{\bf B}^T{\bf S}{\bf A} \label{eqn:L},
\end{align}
where ${\bf{W}}\in {\rm I\!R}^{n\times n}$ and ${\bf{U}}\in {\rm I\!R}^{p\times p}$ are positive definite diagonal weight matrices. $\bf S$ is the solution of the following algebraic Riccati equation,
\begin{equation}
	{\bf S}={\bf A}^T{\bf S}{\bf A}+{\bf W}-{\bf A}^T{\bf S}{\bf B}\left({\bf B}^T{\bf S}{\bf B}+{\bf U}\right)^{-1}{\bf B}^T{\bf S}{\bf A}.
	\label{eqn:S}
\end{equation}
\subsection{Attack Model}
\label{subsec:system}
The schematic diagram of the system under the replay attack is shown in Fig.~\ref{fig:sche_system_attack}. Under the replay attack, the true observations ${\bf y}_k$ are replaced by the delayed version of the observations, \ie, ${\bf z}_k={\bf y}_{k-k_0}$, where $k_0$ represents the delay. The attacker does not need to have any knowledge about the system parameters or the control logic to lunch the replay attack. We assume the attack to start at time $k = \nu$, which is deterministic but unknown. 
\begin{figure}[h!]
	\centering
	\includegraphics[width=55mm]{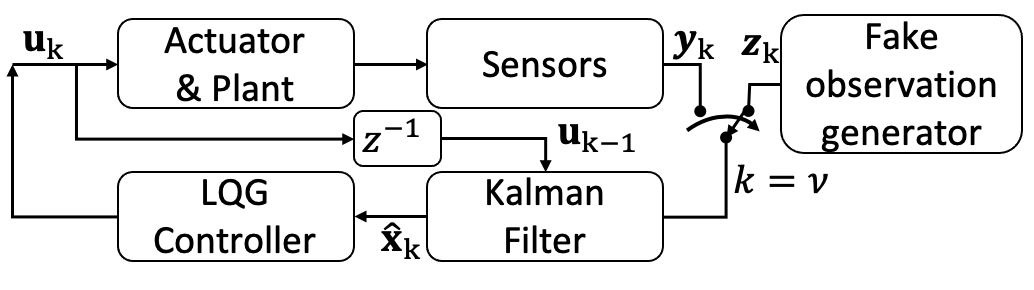}
	\caption{Schematic diagram of the system under replay attack.}
	\label{fig:sche_system_attack.png}
\end{figure}


\section{Replay Attack Detection} \label{sec:cusum}
This section discusses the replay attack detection scheme, the derivations of different parameters and the KLD. It also shows a way to optimize the watermarking signal variance and how to improve the KLD for the systems with relative degree greater than one.
\subsection{Detection Scheme} \label{sec:det_scheme}
To detect the replay attack, we perform the following two main steps. \\
{\bf Step-1}:  Addition of an iid watermarking signal ${\bf e}_k \sim \mathcal{N}(0,{\bf \Sigma}_e)$ to the optimal LQG control input ${\bf u}_k^*$ to generate the control input ${\bf u}_k$ as, 
\begin{equation}
	{\bf u}_k = {\bf u}^*_k+{\bf e}_k.
	\label{eqn:add_ek}
\end{equation}
{\bf Step-2}: Perform the CUSUM test using the joint distributions ${f({{ {\bf \gamma}}_k,{\bf e}_{k-1}})}$ and ${\widetilde f({{\widetilde {\bf \gamma}}_k,{\bf e}_{k-1}})}$ of the innovation signal and the watermarking signal, before and after the attack respectively. We compare the CUSUM statistics with a threshold to select from the following two hypothesises, 
\begin{itemize}
	\item $H_0$: No attack. The estimator receives the true observation ${\bf y}_k$ 
	\item $H_1$: Attack. The estimator receives delayed observation ${\bf z}_k={\bf y}_{k-k_0}$. 
\end{itemize}

The innovation signal $\gamma_k$ during the normal operation of the system is uncorrelated to the watermarking signal, see (\ref{eqn:gamma_normal}). However, on the contrary, the innovation signal $\widetilde\gamma_k$ under the replay attack becomes dependent on the watermarking signal, see (\ref{eqn:gamma_attack}). 
\begin{align}
	\gamma_k&={\bf y}_k-{\bf C}{\bf \hat x}_{k|k-1}  \nonumber \\
	&={\bf C}{\bf A}\left({\bf x}_{k-1}-{\bf \hat x}_{k-1|k-1}\right)+{\bf C}{\bf w}_{k-1}+{\bf v}_k ,
	\label{eqn:gamma_normal} \\
	\widetilde\gamma_k&={\bf y}_{k-k_0}-{\bf C}{\bf \hat x}^F_{k|k-1} \label{eqn:gamma_attack} \\
	&={\bf C}{\bf x}_{k-k_0}+{\bf v}_{k-k_0}-{\bf C}\left({\bf A}+{\bf B}{\bf L}\right){\bf \hat x}^F_{k-1|k-1}-{\bf C}{\bf B}{\bf e}_{k-1}.  \nonumber
\end{align}
The notation ${\bf \hat {x}}^F$ denotes the estimated state when the system is under attack. The CUSUM test statistics $g_k$ is evaluated using Corollary 1.1 from \cite{watermarking_tac} as follows, 
\begin{align}
	&g_k=\max\left(0,g_{k-1}+\log\frac{{\widetilde f}\left({{\bar {\bf \gamma}}_k,{\bf e}_{k-1}}\right)}{f\left({{\bar {\bf \gamma}}_k,{\bf e}_{k-1}}\right)}\right),
	\label{eqn:corr_cusum_subopt} \\
	&\text{where }{\bar \gamma}_k={ \gamma}_k  \text{ before attack, and }  {\bar \gamma}_k={\widetilde \gamma}_k \text{ after attack.} \nonumber \\
	& \text{Furthermore, }  \nonumber \\
	& {\bf \gamma}_{e,k}=\left[{\bf \gamma}_k^T, {\bf e}_{k-1}^T \right]^T  \sim{\cal{N}}\left( {\bf 0},{\bf \Sigma}_{{\gamma}_e} \right),\ { \gamma}_k \sim {\cal{N}}\left( {\bf 0},{\bf \Sigma}_{{\gamma}} \right),  \nonumber \\
	&\text{where } {\bf \Sigma}_{\gamma_e}=\begin{bmatrix}{\bf \Sigma}_\gamma & {\bf 0}_{m\times p} \\{\bf 0}_{p\times m} & {\bf \Sigma}_{e} \end{bmatrix}, \label{eqn:sgima_sq_gamma_e} \text{ and } \\
	& {\bf \widetilde \gamma}_{e,k}=\left[{\bf \widetilde \gamma}_k^T, {\bf e}_{k-1}^T \right]^T \sim{\cal{N}}\left( {\bf 0},{{\bf \Sigma}_{{\widetilde \gamma}_e}}\right),\ {\bf \widetilde \gamma}_k \sim {\cal{N}}\left( {\bf 0},{\bf \Sigma}_{{\widetilde \gamma}} \right), \nonumber \\
	&\text{where } {\bf \Sigma}_{\widetilde \gamma_e}=\begin{bmatrix}{\bf \Sigma}_{\widetilde\gamma} & - {\bf C} {\bf B} {\bf \Sigma}_e \\ - {\bf \Sigma}_e{\bf B}^T {\bf C}^T   & {\bf \Sigma}_{e} \end{bmatrix} \label{eqn:sgima_sq_gamma_e_attack}, \\
	& {\bf \Sigma}_{\gamma}= {\bf C}{\bf P}^T{\bf C}+{\bf R}. \label{eqn:sigma_gamma_sq}
\end{align}
Here, $- {\bf C} {\bf B} {\bf \Sigma}_e$ is the covariance matrix between ${{\widetilde \gamma}_k}$ and ${\bf e}_{k-1}$. The expression of ${\bf \Sigma}_{\widetilde \gamma}$ is provided in Subsection~\ref{subsec:sigma_gamma_sq}.

The threshold for the attack detection is $log(ARL_h)$ where $ARL \le ARL_h$. The decision of attack and no attack is made based on the following,  
\begin{description}
	\item[$H_0:$] Selected, when $g_k < \log(ARL_h)$
	\item[$H_1:$] Selected, when $g_k \ge  \log(ARL_h)$.
\end{description}

The supremum of ADD (SADD), defined as $SADD\triangleq \sup _{1\le\nu<\infty} E_{\nu}\left[ T_{cs}-\nu|T_{cs}>\nu\right]$, will be asymptotically inversely proportion to the KLD, $D\left({\widetilde f},f\right)$, between the two distributions ${\widetilde f}(\cdot)$ and ${f}(\cdot)$ as follows \cite{Tartakovsky2014, Girardin2018}. Here $E_{\nu}(\cdot)$ denotes the expectation with respect to the distribution of the test data when the system is under attack, $T_{cs}$ is the time instant of attack detection, and
\begin{equation}
	SADD
	\rightarrow \frac{\log(ARL_h)}{D\left({\widetilde f},f\right)}, \text{ as } ARL_h \rightarrow \infty.
	\label{ref:subopt_SADD}
\end{equation}
The average run length is defined as, $ARL\triangleq E_{\infty}\left[T_{cs}\right]$, where $E_{\infty}(\cdot)$ denotes the expectation with respect to the distribution of the test data when no attack is present.The KLD, $D\left({\widetilde f},f\right)$, under the CUSUM test will take the following form \cite{watermarking_tac}, 
\begin{align}
	D\left({\widetilde f},f\right) &=\frac{1}{2}\left\{tr\left({\bf \Sigma}_\gamma^{-1} {\bf \Sigma}_{\widetilde \gamma}\right) -m \right. \nonumber \\
	& \left. - \log\frac{\mid{{\bf \Sigma}_{\widetilde \gamma}}-{\bf C}{\bf B}{\bf \Sigma}_e{\bf B}^T{\bf C}^T\mid}{\mid{\bf \Sigma}_{\gamma}\mid}\right\}. \label{eqn:subopt_kld}
\end{align}


\subsection{${\bf \Sigma}_{\widetilde \gamma}$ under Replay Attack} \label{subsec:sigma_gamma_sq}
For the replay attack detection, we assume the attacker's system model to be a partially observed Gauss Markov process (GMP) as follows,
\begin{align}
	{\bf x}_{a,k}&={\bf A}_a{\bf x}_{a,k-1}+{\bf w}_{a,k-1} \label{eqn:hidden_states},  \ \\
	{\bf z}_{k}&={\bf C}_{a}{\bf x}_{a,k} \label{eqn:zk},
\end{align}
where ${\bf x}_{a,k} \in {\rm I\!R}^{n_a}$ and ${\bf w}_{a,k} \sim \mathcal{N}(0,{\bf Q}_a)$ are the hidden state vector and iid noise vector respectively at the $k$-th time instant, and ${\bf Q}_a \in {\rm I\!R}^{n_a \times n_a}$. 
The fake measurements ${\bf z}_k$ will be a delayed version of the true observations ${\bf y}_k$, see \cite{Mo2015}.
Under the replay attack, ${\bf x}_{a,k}$, ${\bf w}_{a,k}$, ${\bf A}_a$, ${\bf C}_a$ and ${\bf Q}_a$ will take the following forms 
\begin{align}
	& {\bf x}_{a,k}=\begin{bmatrix}{\bf x}_{k} & {\hat {\bf x}}_{k|k-1} & {\bf v}_{k} \end{bmatrix}^T \label{eqn:xa_replay} \ \\
	& {\bf w}_{a,k}=\begin{bmatrix}{\bf B}{\bf e}_{k}+{\bf w}_{k} & {\bf B}{\bf e}_{k} & {\bf v}_{k+1} \end{bmatrix}^T \label{eqn:wa_replay} \ \\
	&{\bf A}_a=\cr 
	&\begin{bmatrix}{\bf A}+{\bf B}{\bf L}{\bf K}{\bf C} & {\bf B}{\bf L}\left({\bf I}_n-{\bf K}{\bf C}\right) & {\bf B}{\bf L}{\bf K} \\\left({\bf A}+{\bf B}{\bf L}\right){\bf K}{\bf C} & \left({\bf A}+{\bf B}{\bf L}\right)\left({\bf I}_n-{\bf K}{\bf C}\right) & \left({\bf A}+{\bf B}{\bf L}\right){\bf K} \\{\bf 0} & {\bf 0} & {\bf 0}\end{bmatrix} \label{eqn:Aa_replay} \cr \ \\
	& {\bf C}_a = \begin{bmatrix}{\bf C} & {\bf 0} & {\bf I}_n\end{bmatrix}  \label{eqn:Ca_replay} \ \\
	& {\bf Q}_a=\begin{bmatrix}{\bf B}{\bf \Sigma}_e{\bf B}^T + {\bf Q}& {\bf B}{\bf \Sigma}_e{\bf B}^T & {\bf 0} \\{\bf B}{\bf \Sigma}_e{\bf B}^T & {\bf B}{\bf \Sigma}_e{\bf B}^T & {\bf 0} \\{\bf 0} & {\bf 0} & {\bf R}\end{bmatrix}. \label{eqn:Qa_replay}
\end{align}
${\bf I}_n$ is the identity matrix of size $n$. The derivations follow directly from the system model with watermarking. The ${\bf \Sigma}_{\widetilde \gamma}$, under the replay attack is given as follows, 
\begin{align}
	{\bf \Sigma}_{\widetilde \gamma}&={\bf E}_{zz}(0)-{\bf C}({\bf A}+{\bf B{\bf L)}}{\bf E}_{xz}(-1) \cr
	&-\left[{\bf C}({\bf A}+{\bf B{\bf L)}}{\bf E}_{xz}(-1) \right]^T+{\bf C}{\bf B}{\bf \Sigma}_e{\bf B}^T{\bf C}^T  \cr
	&+{\bf C}({\bf A}+{\bf B}{\bf L}){\bf \Sigma}_{x^Fz}({\bf A}+{\bf B}{\bf L})^T{\bf C}^T \cr
	&+{\bf C}({\bf A}+{\bf B}{\bf L}){\bf \Sigma}_{x^Fe}({\bf A}+{\bf B}{\bf L})^T{\bf C}^T
	\label{eqn:sigma_gamma_attack_replay}
\end{align}
where
\begin{align}
	&{\bf E}_{xz}(-1)= \sum_{i=0}^\infty\mathcal{A}^{i}{\bf K}{\bf C}_a{\bf A}_a^{i+1}{\bf E}_{x_a}\left(0\right){\bf C}_a^T, \label{eqn:Exz1_original_replay}
\end{align}
${\bf E}_{zz}(0) = E\left[{\bf z}_k {\bf z}_k^T\right]$, ${\bf E}_{xz}(-1) = E\left[{\bf x}_{k-1} {\bf z}_k^T\right]$, and ${\bf E}_{x_a}(0) = E\left[{\bf x}_{a,k} {\bf x}_{a,k}^T\right]$. $\mathcal{A} = \left({\bf I}_n-{\bf K}{\bf C}\right)  \left({\bf A}+{\bf B}{\bf L}\right)$. ${\bf \Sigma}_{x^Fz}$ and ${\bf \Sigma}_{x^Fe}$ are the solutions to the following Lyapunov equations,
\begin{align}
	&\mathcal{A}{\bf \Sigma}_{x^Fz}\mathcal{A}^T-{\bf \Sigma}_{x^Fz}+{\bf K}{\bf E}_{zz}(0){\bf K}^T+\mathcal{A}{\bf E}_{xz}(-1){\bf K}^T \nonumber \\
	&+\left(\mathcal{A}{\bf E}_{xz}(-1){\bf K}^T\right)^T = 0 \text{, and} \label{eqn:ExFxF_th1p1_replay} \ \\
	&\mathcal{A}{\bf \Sigma}_{x^Fe}\mathcal{A}^T-{\bf \Sigma}_{x^Fe}+\left({\bf I}_n-{\bf K}{\bf C}\right){\bf B}{\bf \Sigma}_e{\bf B}^T\left({\bf I}_n-{\bf K}{\bf C}\right)^T = 0.
	\label{eqn:ExFxF_th1p2_replay}
\end{align}
Since $\mathcal A$ is assumed to be strictly stable, the Lyapunov equations (\ref{eqn:ExFxF_th1p1_replay}) and (\ref{eqn:ExFxF_th1p2_replay}) will have unique solutions.

\begin{proof}
	The detailed derivation of (\ref{eqn:sigma_gamma_attack_replay}) is provided in Appendix~\ref{apdx:sigma_sq_gamma_attack}.
\end{proof}


\subsection{Optimal watermarking signal variance}
\label{subsec:opt_e}
The addition of watermarking increases the KLD, but at the same time, it also increases the control cost. The increase in LQG control cost, $\Delta LQG$, due to the addition of watermarking is given in \cite{watermarking_tac} as follows,
\begin{equation}
	\Delta LQG=tr\left[\left({\bf B}^T{\bf \Sigma}_L{\bf B}+{\bf U}\right){\bf \Sigma}_e\right], \label{eqn:deltaLQG} 
\end{equation}
where ${\bf \Sigma}_L$ is the solution to the Lyapunov equation
\begin{equation}
	\left({\bf A}+{\bf B}{\bf L}\right)^T{\bf \Sigma}_L\left({\bf A}+{\bf B}{\bf L}\right)-{\bf \Sigma}_L+{\bf L}^T{\bf U}{\bf L} +{\bf W}=0.
	\label{eqn:Sigma_L}
\end{equation} 
Therefore, we want to find the optimal ${\bf \Sigma}_e$ that will maximize the KLD for a given fixed threshold $J$ on the $\Delta LQG$. According to the Theorem~5 from \cite{watermarking_tac}, the optimal ${\bf \Sigma}_e$ will have only one non-zero eigenvalue. Therefore, we search for the optimum ${\bf \Sigma}_e$ within the class of rank one positive semi-definite matrices with the following structure, 
\begin{equation}
	{\bf \Sigma}_e= {\bf v}_{\lambda} {\bf v}_{{\lambda}}^T \text{, where } {\bf v}_{\lambda} = {\sqrt \lambda_e}{\bf v}_e.
	\label{eqn:sigma_e_rank_1_2nd}
\end{equation}
Here, $\lambda_e$ is the non-zero eigenvalue and ${\bf v}_e$ is the corresponding eigenvector. Now, the optimization problem becomes, 
\begin{equation}
	\begin{aligned}
		\max_{{\bf v}_{\lambda}} &\ D\left({\widetilde f},f\right) \\
		\textrm{s.t.}\  & \Delta LQG \le J .
	\end{aligned}
\end{equation}
We have solved the optimization problem using the interior point method \cite{Forsgren2002}. It can also be solved by other non-convex optimizers, such as sequential quadratic programming (SQP) \cite{Boggs1995}, etc. Since the cost function is non-concave, the solution may only be a local optimum.

\subsection{Systems with High Relative Degree}
\label{subsec:relative_degree}
If the system under consideration has a relative degree $d_r = k_e$ where $k_e \ge 2$, then the ${\bf C}{\bf B}{\bf \Sigma}_e{\bf B}^T{\bf C}^T$ term in the KLD expression (\ref{eqn:subopt_kld}) will vanish, which will reduce the overall KLD. In such a situation, the joint distribution of the innovation signal $\gamma_k $ or $\widetilde \gamma_k$ and the delayed version of the watermarking signal, \ie, ${\bf e}_{k-k_e}$, can improve the KLD. Increase in KLD means faster attack detection. The KLD for the joint distribution of the innovation signal and the delayed watermarking signal ${\bf e}_{k-k_e}$ is provided in the following theorem.
\begin{theorem}\label{th:kld_relative_degree}
	If the system has a relative degree of $d_r = k_e$ and the joint distribution of the innovation signal and the watermarking signal ${\bf e}_{k-k_e}$ is considered for the CUSUM test, then the KLD, $D_d\left({\widetilde f},f\right)$, between the normal system and the system under attack will be 
	\begin{align}
		&D_d\left({\widetilde f},f\right)=\frac{1}{2}\left\{tr\left({\bf \Sigma}_\gamma^{-1} {\bf \Sigma}_{\widetilde \gamma}\right) -m \right. \nonumber \\
		&\left.- \log\frac{\mid{{\bf \Sigma}_{\widetilde \gamma}}-{\bf C}{\bf A}^{k_e-1}{\bf B}{\bf \Sigma}_e{\bf B}^T({\bf A}^{k_e-1})^T{\bf C}^T\mid}{\mid{\bf \Sigma}_{\gamma}\mid}\right\}, \label{eqn:kld_relative_degree} 
	\end{align}
\end{theorem}
where the expressions for ${{\bf \Sigma}_{ \gamma}}$ and ${{\bf \Sigma}_{\widetilde \gamma}}$ will be the same as in (\ref{eqn:sigma_gamma_sq}) and (\ref{eqn:sigma_gamma_attack_replay}) respectively but $E\left[\gamma_k{\bf e}_{k-k_e}^T\right]$ and $E\left[\widetilde \gamma_k{\bf e}_{k-k_e}^T\right]$ need to be derived. However, the expression of ${{\bf \Sigma}_{\widetilde \gamma}}$ can be simplified using the information ${\bf C}{\bf B}=0$ for a system with $d_r \ge 2$. 

\begin{proof}
	The proof of Theorem~\ref{th:kld_relative_degree} is provided in Appendix~\ref{apdx:kld_relative_degree_2}.
\end{proof}

\section{Numerical Results}
\label{sec:numerical_results}
In this section, we illustrate the replay attack detection methodology proposed in this paper using three different system models. The three different systems are System-A: A second-order open-loop unstable multiple inputs and single output (MISO) system, System-B: A fourth-order open-loop stable multiple inputs and multiple outputs (MIMO) system, and System-C: A second-order open-loop unstable MISO system with relative degree two. The system parameters are provided in Appendix~\ref{apdx:system_param}. System-B is a linearized minimum phase quadruple tank system taken from \cite{Johansson1998}. Only the level sensor gains are increased to make the magnitude of the product ${\bf C}{\bf B}$ numerically significant. 
\subsection{Replay attack detection}
\begin{figure}[h!]
	\centering
	\includegraphics[width=\figwidth]{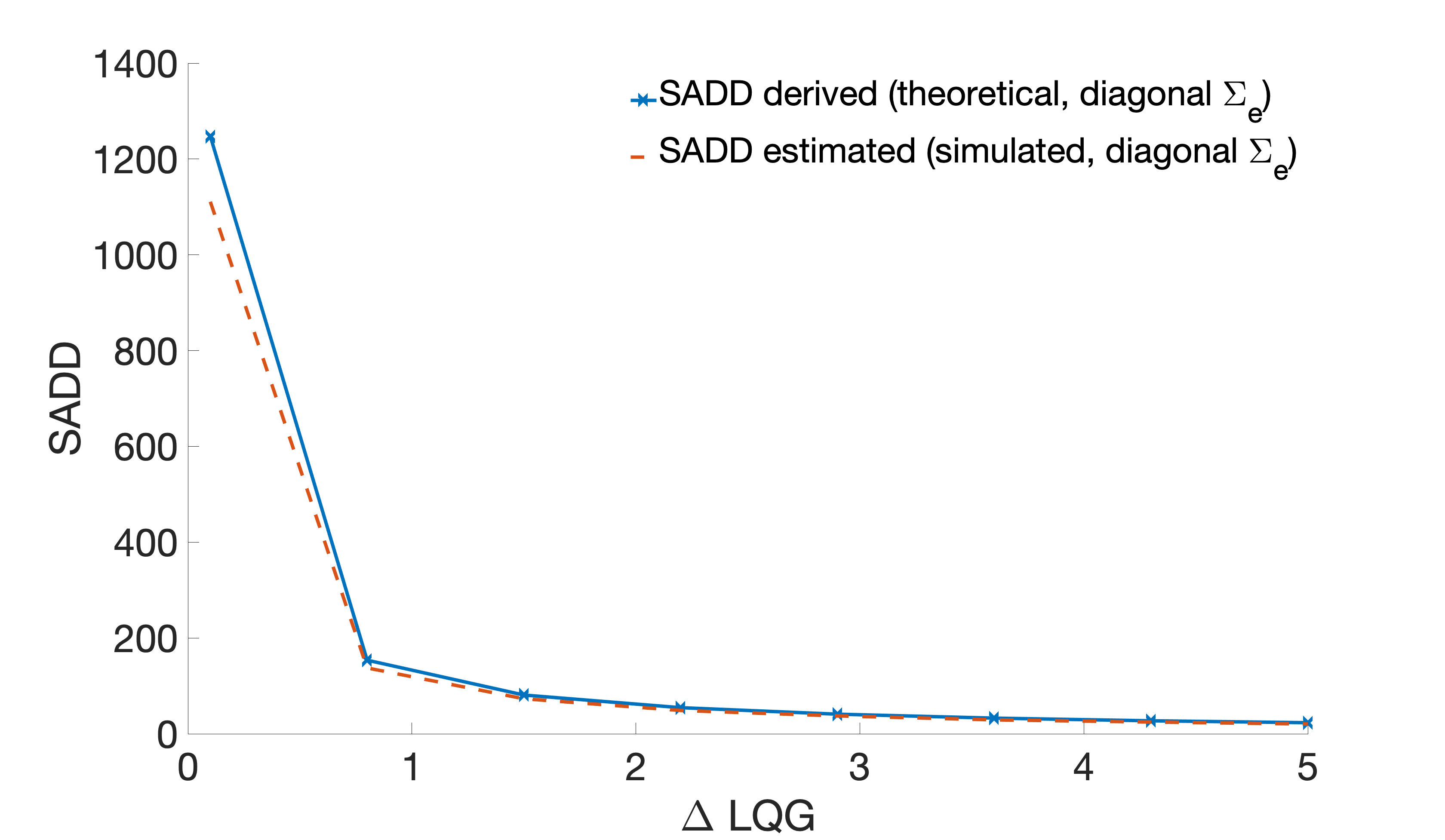}
	\caption{SADD vs. $\Delta LQG$ plot for System-B under replay attack.}
	\label{fig:replay_attack_mimo}
\end{figure}

Figure~\ref{fig:replay_attack_mimo} shows the tradeoff between the SADD and the increase in $\Delta LQG$ when the System-B is under a replay attack. We plot the derived SADD using the theory developed in this paper, and the estimated SADD using the simulated data where ${\bf \Sigma}_e$ is assumed to be diagonal and all the watermarking signals have equal power. Therefore, we can claim that the proposed sequential detection technique can detect replay attacks. Figure~\ref{fig:replay_attack_mimo} also illustrates that it is hard to detect a replay attack with low watermarking signal power. This is implicit in SADD, but there is a sharp increase before $\Delta LQG \approx 0.8$, which corresponds to $\Sigma_e = \text{diag}\begin{bmatrix}0.29 & 0.29 \end{bmatrix}$, is observed for the system model under consideration. 
\subsection{Optimum an  non-optimum ${\bf \Sigma}_e$}
\begin{figure}[h!]
	\centering
	\includegraphics[width=\figwidth]{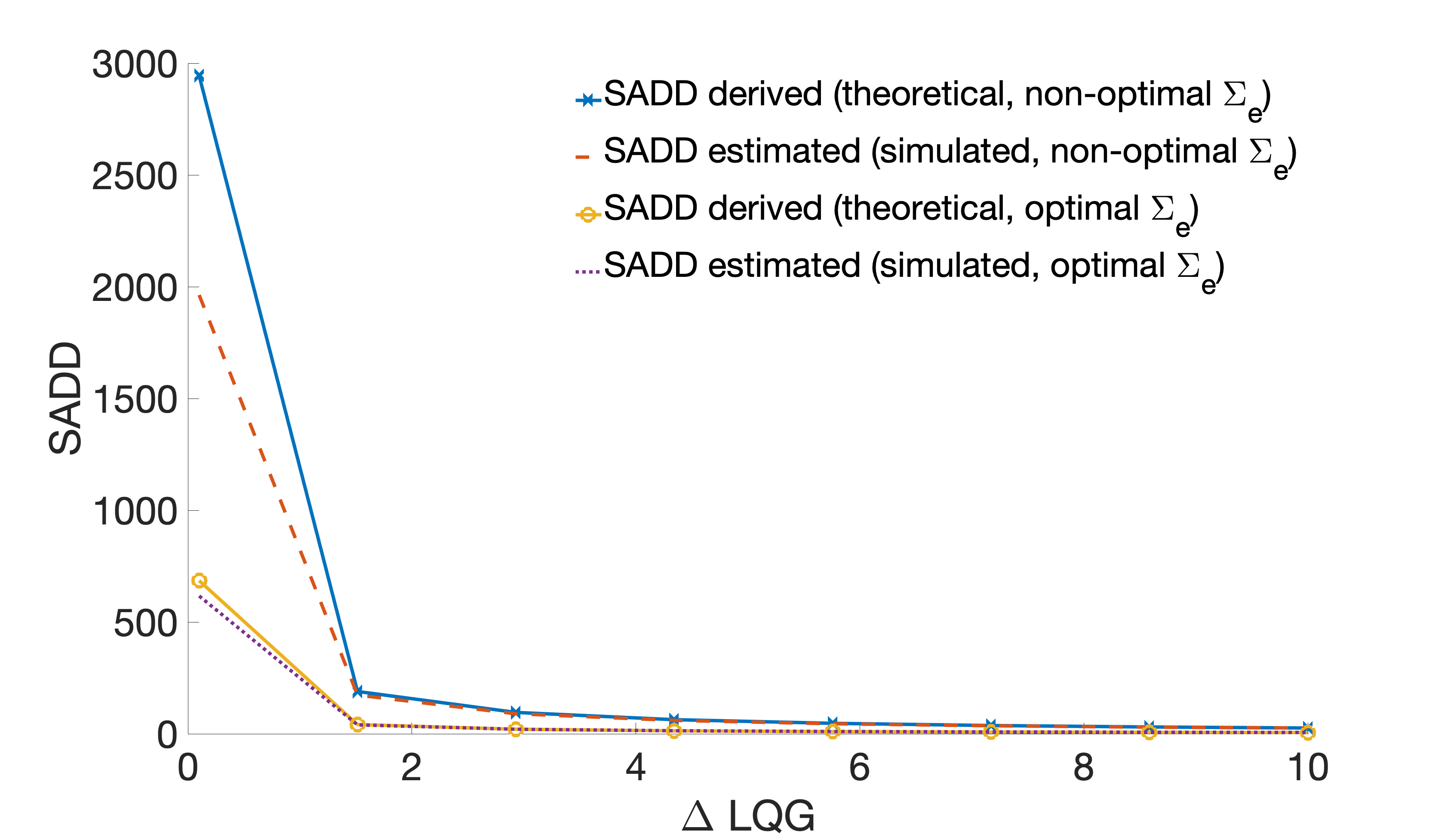}
	\caption{SADD vs. $\Delta LQG$ plot for System-A under replay attack with optimum and non-optimum ${\bf \Sigma}_e$.}
	\label{fig:replay_attack_miso}
\end{figure}

Figure~\ref{fig:replay_attack_miso} shows the SADD vs $\Delta LQG$ plots for System-A using the optimized ${\bf \Sigma}_e$ and a diagonal ${\bf \Sigma}_e$ with equal signal power when the system is under replay attack. It is evident that optimizing ${\bf \Sigma}_e$ improves SADD for a fixed upper threshold on $\Delta LQG$. We can also say that the same  SADD can be achieved for much reduced $\Delta LQG$.
\subsection{System with higher relative degree}
\begin{figure}[h!]
	\centering
	\includegraphics[width=\figwidth]{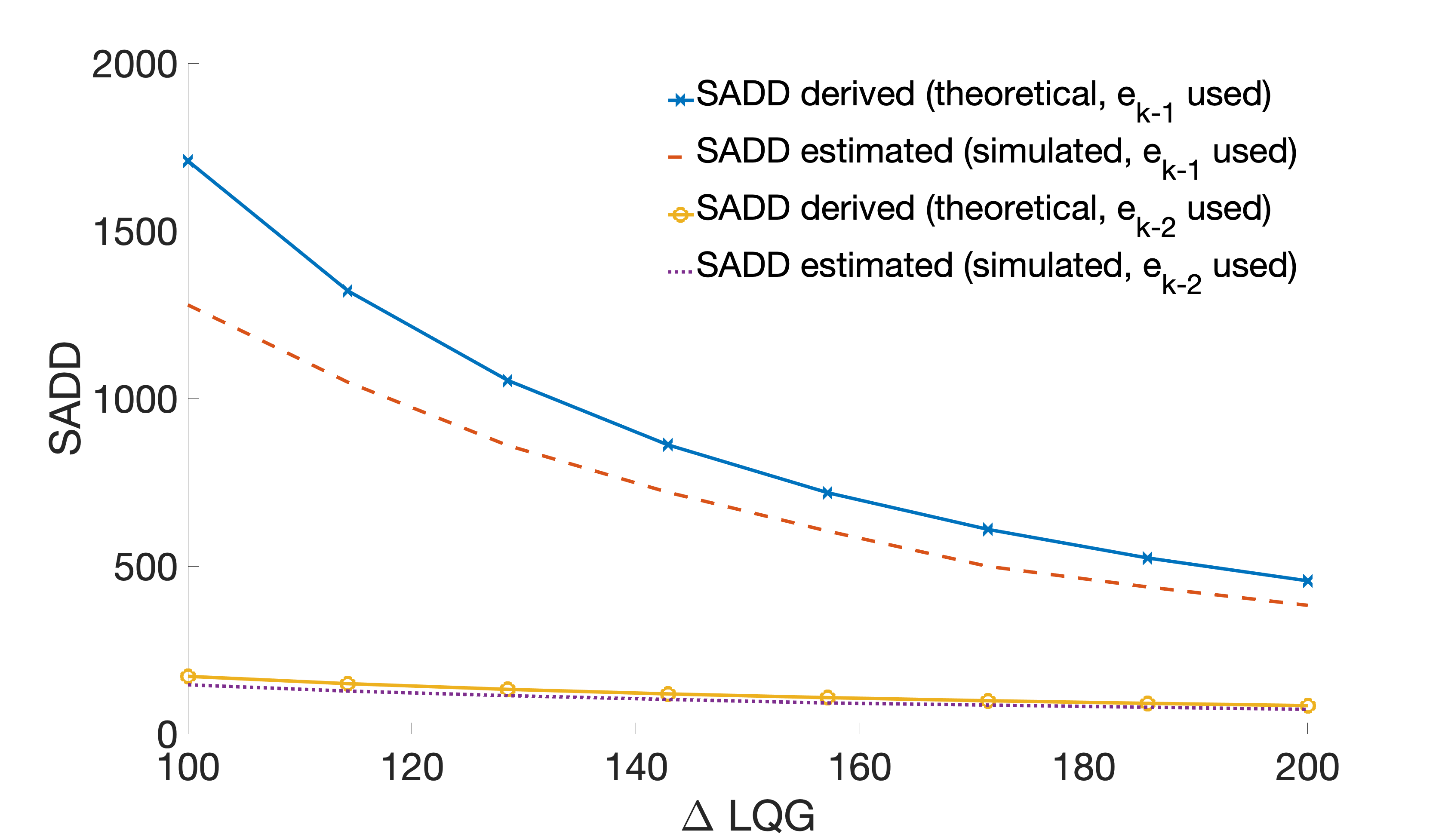}
	\caption{SADD vs. $\Delta LQG$ plots for System-C with relative degree 2.}
	\label{fig:add_vs_Dlqg_system_rd2}
\end{figure}

Figure~\ref{fig:add_vs_Dlqg_system_rd2} shows the benefit of using the delayed version of watermarking signal, \ie, ${\bf e}_{k-k_e}$ for a system with relative degree $d_r = k_e$ as discussed in Theorem~\ref{th:kld_relative_degree}. System-C with relative degree $d_r = 2$ is used to generate the plots of Fig.~\ref{fig:add_vs_Dlqg_system_rd2}. We can see reductions in $\Delta LQG$ to achieve the same $SADD$ between any two points on the $\Delta LQG$ axis.
\subsection{Comparison with optimal NP detector}
\begin{figure}[h!]
	\centering
	\includegraphics[width=\figwidth]{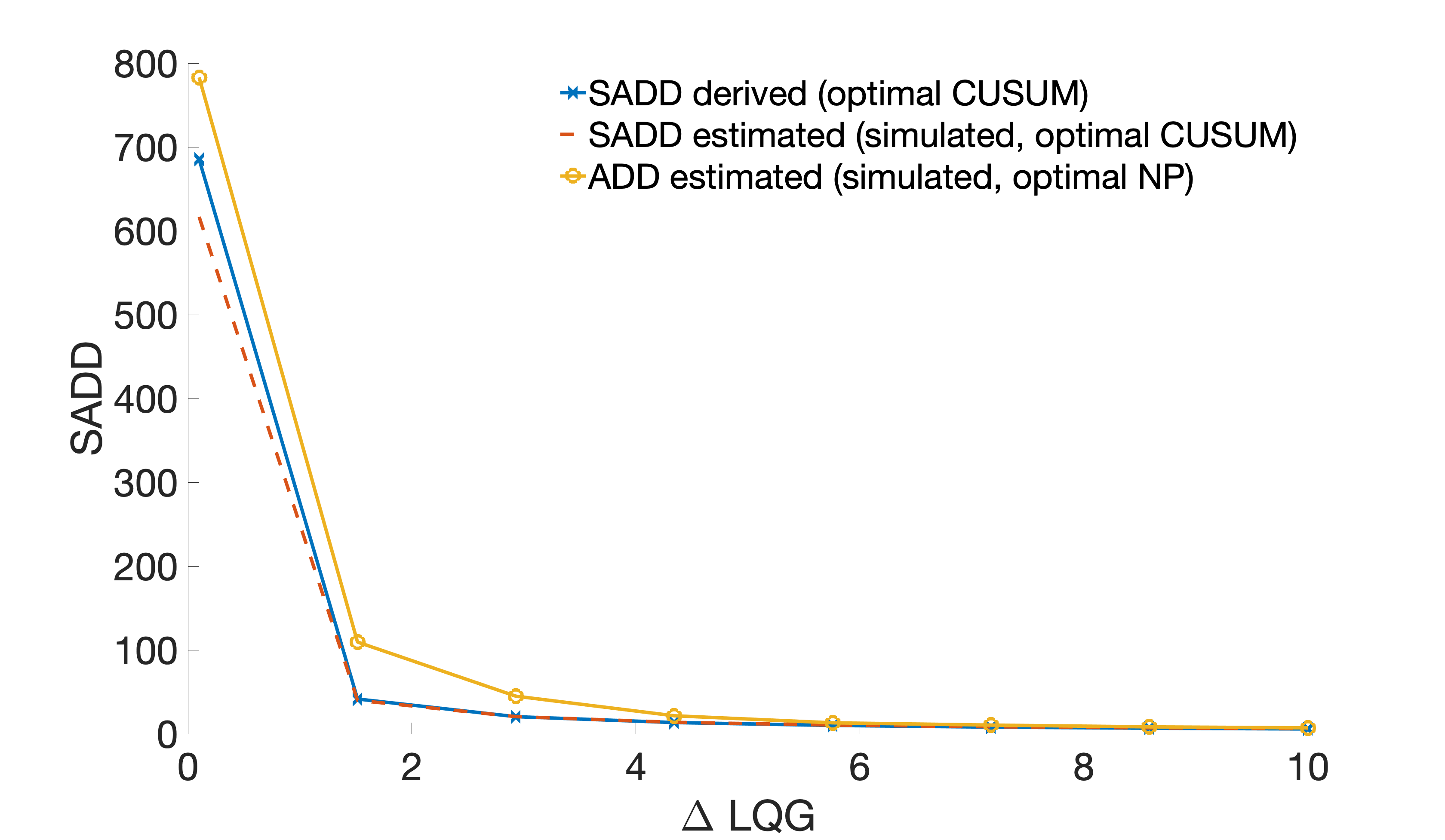}
	\caption{SADD vs. $\Delta LQG$ plot for System-A under CUSUM test and NP test.}
	\label{fig:compare_cusum_vs_chi_sq}
\end{figure}

Figure~\ref{fig:compare_cusum_vs_chi_sq} shows the tradeoff between the ADD and the increase in $\Delta LQG$ for System-A under the CUSUM test and the method reported in \cite{Mo2015} for the detection of replay attacks. We plot the derived SADD using the theory developed in this paper, the estimated SADD applying the CUSUM test on the simulated data, and the estimated ADD applying the test reported in \cite{Mo2015} on the simulated data. It is clear from the figure that we can achieve lower detection delay for the same LQG loss with the method proposed in this paper compared to the one reported in \cite{Mo2015}. FAR is taken to be the reciprocal of ARL \cite{Giraldo2019}.     



\section{Conclusion}
\label{sec:conclusion}
We have addressed the problem of resilient replay attack detection using the CUSUM test. The detection delay and the corresponding loss in the control cost are studied. The KLD expression between the distributions before the attack and after the replay attack is derived. SADD is asymptotically inversely proportional to the KLD. The KLD reduces for the systems with relative-degree higher than one. We have shown a technique of using a delayed version of the watermarking to improve the KLD for such systems. The results shown are in close agreement with the theory presented in the paper. We have also discussed a way to optimize the watermarking signal variance to maximize the KLD under the replay attack for a fixed increase in the control cost.
\appendices

\section{Derivation of  ${\bf \Sigma}_{\widetilde \gamma}$}
\label{apdx:sigma_sq_gamma_attack}
It is assumed that the fake observations will be generated by the following partially observed GMP
\begin{align}
	{\bf x}_{a,k}&={\bf A}_a{\bf x}_{a,k-1}+{\bf w}_{a,k-1} \label{eqn:hidden_states_apndx} \ \\
	{\bf z}_{k}&={\bf C}_{a}{\bf x}_{a,k} \label{eqn:zk_apndx}
\end{align}
where ${\bf x}_{a,k} \in {\rm I\!R}^{n_a}$ and ${\bf w}_{a,k} \sim \mathcal{N}(0,{\bf Q}_a)$ are the hidden state vector and iid noise vector respectively at the $k$-th time instant, and ${\bf Q}_a \in {\rm I\!R}^{n_a \times n_a}$.

Since the true measurement ${\bf y}_k$ is stationary, the attacker will keep the fake measurement ${\bf z}_k$ stationary by taking the initial state covariance as ${\bf E}_{x_a}(0)$ to remain stealthy, where ${\bf E}_{x_a}(0)$ is the solution of the following Lyapunov equation. 
\begin{align}
	{\bf E}_{x_a}(0)={\bf A}_a{\bf E}_{x_a}(0){\bf A}^T_a+{\bf Q}_a.
	\label{eqn:cxo_apndx}
\end{align}

The variance of the innovation signal (${\bf \Sigma}_{\widetilde \gamma}$) when the system is under attack is derived as follows. 
\begin{align}
	{\bf \Sigma}_{\widetilde \gamma}=E\left[{\widetilde \gamma}_k{\widetilde \gamma}_k^T\right]. 
\end{align}
Using (\ref{eqn:gamma_attack}), and applying the knowledge that ${\bf e}_{k-1}$ is uncorrelated with ${\bf z}_k$ and ${\hat {\bf x}}_{k-1|k-1}^F$, we get the following expression of ${\bf \Sigma}_{\widetilde \gamma}$,
\begin{align}
	&{\bf \Sigma}_{\widetilde \gamma} = E\left[{\bf z}_k{\bf z}_k^T\right]-{\bf C}\left({\bf A}+{\bf B}{\bf L}\right)E\left[{\hat {\bf x}}_{k-1|k-1}^F{\bf z}_k^T\right] \cr
	&-\left({\bf C}\left({\bf A}+{\bf B}{\bf L}\right)E\left[{\hat {\bf x}}_{k-1|k-1}^F{\bf z}_k^T\right]\right)^T+{\bf C}{\bf B}{\bf \Sigma}_e{\bf B}^T{\bf C}^T \cr
	&+{\bf C}\left({\bf A}+{\bf B}{\bf L}\right)E\left[{\hat {\bf x}}_{k-1|k-1}^F\left({\hat {\bf x}}_{k-1|k-1}^{F}\right)^T\right]\left({\bf A}+{\bf B}{\bf L}\right)^T{\bf C}^T. \cr
	\label{eqn:Sigma_gamma_attack_stage_1}
\end{align}
$E\left[{\hat {\bf x}}_{k-1|k-1}^F{\bf z}_k^T\right]$ is calculated as follows, 
\begin{align}
	{\hat {\bf x}}_{k-1|k-1}^F&={\hat {\bf x}}_{k-1|k-2}^F+{\bf K}\left({\bf z}_{k-1} -{\bf C}{\hat {\bf x}}_{k-1|k-2}^F\right) \cr
	&={\bf K}{\bf z}_{k-1}+\mathcal{A}{\hat {\bf x}}_{k-2|k-2}^F \cr
	&+\left({\bf I}_n-{\bf K}{\bf C}\right){\bf B}{\bf e}_{k-2},\text{ [using  (\ref{eqn:add_ek}]}, \cr
	\text{where }\mathcal{A}&=\left({\bf I}_n-{\bf K}{\bf C}\right)\left({\bf A}+{\bf B}{\bf L}\right).
	\label{eqn:xF_k_1_k_1}
\end{align}
We define ${\bf E}_{xz}\left(-k_0\right) \triangleq E\left[{\hat {\bf x}}_{k-k_0|k-k_0}^F{\bf z}_k^T\right]$, 
\begin{align}
	\begin{split}
		&=E\left[\left({\bf K}{\bf z}_{k-k_0}+\mathcal{A}{\hat {\bf x}}_{k-k_0-1|k-k_0-1}^F \right.\right.  \\
		&+\left.\left.\left({\bf I}_n-{\bf K}{\bf C}\right){\bf B}{\bf e}_{k-k_0-1}\right){\bf z}_{k}^T \right] \text{, [using (\ref{eqn:xF_k_1_k_1})]}  \cr
		&={\bf K}{\bf E}_{zz}\left(-k_0\right)+\mathcal{A}{\bf E}_{xz}\left(-k_0-1\right),
	\end{split}
	\label{eqn:Exz_k0}
\end{align}
where ${\bf e}_{k-k_0-1}$ and ${\bf z}_{k}$ are uncorrelated, and ${\bf E}_{zz}\left(-k_0\right)=E\left[{\bf z}_{k-k_0}{\bf z}_{k}^T\right]$. ${\bf E}_{x_a}\left(-k_0\right)$ is the correlation between the states of attack system, ${\bf x}_{a,k-k_0}$ and ${\bf x}_{a,k}$, which is evaluated as follows. 
\begin{align}
	&{\bf E}_{x_a}\left(-k_0\right)={\bf E}_{x_a}\left(k_0\right)=E\left[{\bf x}_{a,k}{\bf x}_{a,k-k_0}^T\right], \cr
	&{\bf E}_{x_a}\left(-1\right)=E\left[{\bf A}_a{\bf x}_{a,k-1}{\bf x}_{a,k-1}^T+{\bf w}_{a,k-1}{\bf x}_{a,k-1}^T\right] \cr
	&{\bf E}_{x_a}\left(-1\right)={\bf A}_a{\bf E}_{x_a}\left(0\right) \text{, [}{\bf w}_{a,k} \text{ and }{\bf x}_{a,k} \text{ uncorrelated]}.\cr
	&\text{Similarly,} \cr
	&{\bf E}_{x_a}\left(-2\right)={\bf A}_a{\bf E}_{x_a}\left(-1\right)={\bf A}_a^2{\bf E}_{x_a}\left(0\right) \text{, and} \cr
	&{\bf E}_{x_a}\left(-k_0\right) = {\bf A}_a^{k_0}{\bf E}_{x_a}\left(0\right).
	\label{eqn:Exa_k0}
\end{align}
The system matrix ${\bf A}_a$ is assumed to be stable because the attacker will always try to generate fake observations which are bounded and will mimic the true observations to remain stealthy. For a stable ${\bf A}_a$, 
\begin{align}
	&{\bf A}_a^{k_0} \rightarrow 0 \text{, as } k_0 \rightarrow \infty \cr
	&\text{therefore, }{\bf E}_{x_a}\left(-k_0\right) \rightarrow 0 \text{, as } k_0 \rightarrow \infty.
	\label{eqn:Exa_infty}
\end{align}
The expression for ${\bf E}_{zz}(-k_0)={\bf E}_{zz}(k_0)=E\left[{\bf z}_{k-k_0}{\bf z}_{k}^T\right]$ is derived as 
\begin{align}
	{\bf E}_{zz}(-k_0)&={\bf C}_a{\bf E}_{x_a}(-k_0){\bf C}_a^T \text{, [using (\ref{eqn:zk})]} \cr
	&={\bf C}_a{\bf A}_a^{k_0}{\bf E}_{x_a}\left(0\right){\bf C}_a^T \text{, [using (\ref{eqn:Exa_infty})]}.
	\label{eqn:Ezz_k0}
\end{align}
Using (\ref{eqn:Exz_k0}) and (\ref{eqn:Ezz_k0}), we can write the expression of ${\bf E}_{xz}(-1)$ as
\begin{align}
	{\bf E}_{xz}&\left(-1\right)={\bf K}{\bf E}_{zz}\left(-1\right)+\mathcal{A}{\bf E}_{xz}\left(-2\right) \cr
	&={\bf K}{\bf C}_a{\bf A}_a^{}{\bf E}_{x_a}\left(0\right){\bf C}_a^T+\mathcal{A}\left({\bf K}{\bf E}_{zz}\left(-2\right)+\mathcal{A}{\bf E}_{xz}\left(-3\right) \right) \cr  
	& \text{, [by replacing } {\bf E}_{xz}\left(-2\right) \text{ using (\ref{eqn:Exz_k0})]} \cr
	&={\bf K}{\bf C}_a{\bf A}_a^{}{\bf E}_{x_a}\left(0\right){\bf C}_a^T+\mathcal{A}{\bf K}{\bf C}_a{\bf A}_a^{2}{\bf E}_{x_a}\left(0\right){\bf C}_a^T \cr 
	&+\mathcal{A}^2{\bf E}_{xz}\left(-3\right).
	\label{eqn:Exz_1}
\end{align}
Repeating the same technique, ${\bf E}_{xz}\left(-1\right)$ will take the following form,
\begin{align}
	{\bf E}_{xz}\left(-1\right)
	&={\bf K}{\bf C}_a{\bf A}_a^{}{\bf E}_{x_a}\left(0\right){\bf C}_a^T+\mathcal{A}{\bf K}{\bf C}_a{\bf A}_a^{2}{\bf E}_{x_a}\left(0\right){\bf C}_a^T \cr
	&+\mathcal{A}^2{\bf K}{\bf C}_a{\bf A}_a^{3}{\bf E}_{x_a}\left(0\right){\bf C}_a^T+\cdots \cr
	&\sum_{i=0}^\infty\mathcal{A}^{i}{\bf K}{\bf C}_a{\bf A}_a^{i+1}{\bf E}_{x_a}\left(0\right){\bf C}_a^T.
	\label{eqn:Exz_1_series}
\end{align}
${\bf E}_{xz}\left(-1\right)$ can be evaluated numerically by taking a large number of terms for the summation (\ref{eqn:Exz_1_series}), until the rest of the terms become negligible. 
$E\left[{\hat {\bf x}}_{k-1|k-1}^F\left({\hat {\bf x}}_{k-1|k-1}^{F}\right)^T\right]$ is evaluated as follows using (\ref{eqn:xF_k_1_k_1}).
\begin{align}
	&{\bf E}_{x^Fx^F}(0)=E\left[{\hat {\bf x}}_{k-1|k-1}^F\left({\hat {\bf x}}_{k-1|k-1}^{F}\right)^T\right] \cr
	&={\bf K}E\left[{\bf z}_{k-1}{\bf z}_{k-1}^T\right]{\bf K}^T+\mathcal{A}E\left[{\hat {\bf x}}_{k-2|k-2}^F{\bf z}_{k-1}^T\right]{\bf K}^T \cr
	&+\left(\mathcal{A}E\left[{\hat {\bf x}}_{k-2|k-2}^F{\bf z}_{k-1}^T\right]{\bf K}^T\right)^T  \cr
	& +\mathcal{A}E\left[{\hat {\bf x}}_{k-2|k-2}^F\left({\hat {\bf x}}_{k-2|k-2}^{F}\right)^T\right]\mathcal{A}^T \cr
	&+\left({\bf I}_n-{\bf K}{\bf C}\right){\bf B}E\left[{\bf e}_{k-2}{\bf e}_{k-2}^T\right]{\bf B}^T\left({\bf I}_n-{\bf K}{\bf C}\right)^T.
	\label{eqn:ExFxF}
\end{align}
Therefore, ${\bf E}_{x^Fx^F}(0)$ is the solution to the following Lyapunov equation,
\begin{align}
	&\mathcal{A}{\bf E}_{x^Fx^F}(0)\mathcal{A}^T-{\bf E}_{x^Fx^F}(0)+{\bf K}{\bf E}_{zz}(0){\bf K}^T \cr
	&+\mathcal{A}{\bf E}_{xz}(-1){\bf K}^T+\left(\mathcal{A}{\bf E}_{xz}(-1){\bf K}^T\right)^T \cr
	&+\left({\bf I}_n-{\bf K}{\bf C}\right){\bf B}{\bf \Sigma}_e{\bf B}^T\left({\bf I}_n-{\bf K}{\bf C}\right)^T = 0 \text{, [(\ref{eqn:Ezz_k0}) used]}. \cr
\end{align}
${\bf E}_{x^Fx^F}(0)$ is divided into two parts, ${\bf \Sigma}_{x^Fz}$ and ${\bf \Sigma}_{x^Fe}$ which are independent of the watermarking signal and the fake observations, respectively. ${\bf \Sigma}_{x^Fz}$ and ${\bf \Sigma}_{x^Fe}$ are the solution to the following Lyapunov equations, 
\begin{align}
	&\mathcal{A}{\bf \Sigma}_{x^Fz}\mathcal{A}^T-{\bf \Sigma}_{x^Fz}+{\bf K}{\bf E}_{zz}(0){\bf K}^T+\mathcal{A}{\bf E}_{xz}(-1){\bf K}^T \cr
	&+\left(\mathcal{A}{\bf E}_{xz}(-1){\bf K}^T\right)^T = 0 \text{,} \cr
	&\mathcal{A}{\bf \Sigma}_{x^Fe}\mathcal{A}^T-{\bf \Sigma}_{x^Fe}+\left({\bf I}_n-{\bf K}{\bf C}\right){\bf B}{\bf \Sigma}_e{\bf B}^T\left({\bf I}_n-{\bf K}{\bf C}\right)^T = 0 \text{,} \cr
	&\text{and }{\bf E}_{x^Fx^F}(0)={\bf \Sigma}_{x^Fz}+{\bf \Sigma}_{x^Fe}. 
	\label{eqn:ExFxF_apn}
\end{align}
Using (\ref{eqn:Ezz_k0}) and (\ref{eqn:ExFxF_apn}), we can rewrite the expression for ${\bf \Sigma}_{\widetilde \gamma} $ as,
\begin{align}
	{\bf \Sigma}_{\widetilde \gamma}&={\bf E}_{zz}(0)-{\bf C}({\bf A}+{\bf B{\bf L)}}{\bf E}_{xz}(-1) \cr
	&-\left[{\bf C}({\bf A}+{\bf B{\bf L)}}{\bf E}_{xz}(-1) \right]^T+{\bf C}{\bf B}{\bf \Sigma}_e{\bf B}^T{\bf C}^T  \cr
	&+{\bf C}({\bf A}+{\bf B}{\bf L}){\bf \Sigma}_{x^Fz}({\bf A}+{\bf B}{\bf L})^T{\bf C}^T \cr
	&+{\bf C}({\bf A}+{\bf B}{\bf L}){\bf \Sigma}_{x^Fe}({\bf A}+{\bf B}{\bf L})^T{\bf C}^T.
\end{align}


\section{Proof of Theorem~\ref{th:kld_relative_degree}}
\label{apdx:kld_relative_degree_2}
Assumed: System has relative degree $d_r = k_e$, \ie,
\begin{equation}
	{\bf C}{\bf A}^{i}{\bf B} = {\bf 0} \text{, for } i < k_e-1.
	\label{eqn:cabeq0}
\end{equation}
Since $\gamma_k$ before the attack is iid, $E\left[\gamma_k{\bf e}_{k-k_e}^T\right]={\bf 0}$. 
Applying (\ref{eqn:cabeq0}) in (\ref{eqn:gamma_attack}), and denoting ${\bf z}_k = {\bf y}_{k - k_0}$, we get,
\begin{equation}
	\widetilde \gamma_k={\bf z}_k-{\bf C}{\bf A}{\bf \hat x}_{k-1|k-1}^F. 
	\label{eqn:gamma_attack_rd}
\end{equation}
Therefore, 
\begin{align}
	&E\left[\widetilde \gamma_k {\bf e}_{k-k_e}^T \right]=-{\bf C}{\bf A}E\left[{\bf \hat x}_{k-1|k-1}^F{\bf e}_{k-k_e}^T \right], \text{ and}  \label{eqn:Egammaek0} \\\ 
	&{\bf \hat x}_{k-1|k-1}^F = {\bf K}{\bf z}_{k-1}+\mathcal{A}{\bf \hat x}_{k-2|k-2}^F+{\bf B}{\bf e}_{k-2} \label{eqn:xhfkm1},
\end{align}
where ${\bf e}_{k-k_e}$ is uncorrelated to ${\bf z}_k$. Using (\ref{eqn:xhfkm1}) recursively, we derive
\begin{equation}
	\begin{aligned}
		&{\bf \hat x}_{k-1|k-1}^F = \sum_{i=2}^{k_e}{\mathcal {A}^{i-2}}{\bf K}{\bf \bar y}_{k-i+1}+\mathcal{A}^{k_e-1}{\bf \hat x}_{k-k_e|k-k_e}^F  \\
		&+\sum_{i=2}^{k_e}{\mathcal {A}^{i-2}}{\bf B}{\bf e}_{k-i}, \\
		&\text{where } {\bf \bar y}_k= {\bf y}_k \text{ if }k<\nu\text{, and } {\bf \bar y}_k ={\bf z}_k\text{ otherwise.}
		\label{eqn:xhf_k0}
	\end{aligned}	
\end{equation}
Applying (\ref{eqn:xhf_k0}) in (\ref{eqn:Egammaek0}) and using (\ref{eqn:cabeq0}), we get
\begin{align}
	E\left[\widetilde \gamma_k {\bf e}_{k-k_e}^T \right]=-{\bf C}{\bf A}{\mathcal {A}^{k_e-2}}{\bf B}{\bf \Sigma}_e \label{eqn:eqn:Egammaek0_2}
\end{align}
where ${\bf e}_{k-k_e}$ is uncorrelated to ${\bf \hat x}_{k-k_e|k-k_e}^F$ and ${\bf z}_{k-i+1}$. Applying multinomial theorem on $\mathcal{A}^{k_e-1} = \left({\bf A}-{\bf K}{\bf C}{\bf A}+{\bf B}{\bf L} \right)^{k_e-1}$ and using (\ref{eqn:cabeq0}), we get 
\begin{align}
	E\left[\widetilde \gamma_k {\bf e}_{k-k_e}^T \right]=-{\bf C}{\bf A}^{k_e-1}{\bf B}{\bf \Sigma}_e \label{eqn:eqn:Egammaek0_3}.
\end{align}

\section{System Parameters}
\label{apdx:system_param}
For the System-A and System-B, $ARL_h = 1000$. For the System-C, $ARL_h = 100$\\
\textbf{System-A parameters}:
\begin{align*}
	{\bf A} &=\begin{bmatrix}0.75 & 0.2 \\0.2 & 1.0 \end{bmatrix}           &  {\bf B} &=\begin{bmatrix}0.9 & 0.5 \\0.1 & 1.2 \end{bmatrix}              &  {\bf C}&=\begin{bmatrix}1.0 & -1.0  \end{bmatrix} \\
	{\bf Q} &=diag\begin{bmatrix}1 & 1  \end{bmatrix}           &  {\bf R} &=1              &  {\bf W}&=diag\begin{bmatrix}1 & 2  \end{bmatrix} \\
	{\bf U} &=diag\begin{bmatrix}0.4 & 0.7  \end{bmatrix}           &  \sigma_z^2 &=10             &  \rho&=0.5 
\end{align*}
\textbf{System-B parameters}:
\begin{align*}
	{\bf A} =\begin{bmatrix}0.9683 &0&0.0819 &0 \\ 0&0.9780&0&0.06377 \\ 0&0&0.9167&0 \\ 0&0&0&0.9355 \end{bmatrix} 
\end{align*}
\begin{align*}
	{\bf B} &=\begin{bmatrix} 0.1638&0.004 \\0.002&0.1242 \\ 0&0.0917 \\ 0.0604&0  \end{bmatrix}  & {\bf C} &= \begin{bmatrix} 5 &0 &0 &0 \\  0 &5 &0 &0  \end{bmatrix} \\
	{\bf Q} &=diag\begin{bmatrix}0.25 & 0.25 & 0.25 & 0.25  \end{bmatrix}           &  {\bf R} &=diag\begin{bmatrix}0.5 & 0.5  \end{bmatrix} \\
	{\bf W} &=diag\begin{bmatrix}5 & 5 & 1 & 1  \end{bmatrix}           &  {\bf U} &=diag\begin{bmatrix}2 & 2  \end{bmatrix} \\
	{\bf A}_a &=diag\begin{bmatrix}0.4 & 0.2 & 0.2 & 0.7  \end{bmatrix}           &  {\bf Q}_a &=diag\begin{bmatrix}5 & 5  \end{bmatrix} 
\end{align*}
\textbf{System-C parameters}:
\begin{align*}
	{\bf B} &=\begin{bmatrix}0.9 & 0.5 \\1.3 & 0.72 \end{bmatrix}              &  {\bf C}&=\begin{bmatrix}1.3 & -0.9  \end{bmatrix} 
\end{align*}
The rest of the parameters are same as System-A.

%
%

\bibliographystyle{IEEEtran}
\bibliography{IEEEabrv,Postdoc}

\begin{thebibliography}{10}
\providecommand{\url}[1]{#1}
\csname url@samestyle\endcsname
\providecommand{\newblock}{\relax}
\providecommand{\bibinfo}[2]{#2}
\providecommand{\BIBentrySTDinterwordspacing}{\spaceskip=0pt\relax}
\providecommand{\BIBentryALTinterwordstretchfactor}{4}
\providecommand{\BIBentryALTinterwordspacing}{\spaceskip=\fontdimen2\font plus
\BIBentryALTinterwordstretchfactor\fontdimen3\font minus
  \fontdimen4\font\relax}
\providecommand{\BIBforeignlanguage}[2]{{%
\expandafter\ifx\csname l@#1\endcsname\relax
\typeout{** WARNING: IEEEtran.bst: No hyphenation pattern has been}%
\typeout{** loaded for the language `#1'. Using the pattern for}%
\typeout{** the default language instead.}%
\else
\language=\csname l@#1\endcsname
\fi
#2}}
\providecommand{\BIBdecl}{\relax}
\BIBdecl

\bibitem{Alguliyev2018}
\BIBentryALTinterwordspacing
R.~Alguliyev, Y.~Imamverdiyev, and L.~Sukhostat, ``{Cyber-physical systems and
  their security issues},'' \emph{Comput. Ind.}, vol. 100, no. July 2017, pp.
  212--223, 2018. [Online]. Available:
  \url{https://doi.org/10.1016/j.compind.2018.04.017}
\BIBentrySTDinterwordspacing

\bibitem{Langner2011}
R.~Langner, ``{Stuxnet: Dissecting a cyberwarfare weapon},'' \emph{IEEE Secur.
  Priv.}, vol.~9, no.~3, pp. 49--51, 2011.

\bibitem{Mo2015}
Y.~Mo, S.~Weerakkody, and B.~Sinopoli, ``{Physical authentication of control
  systems: Designing watermarked control inputs to detect counterfeit sensor
  outputs},'' \emph{IEEE Control Syst.}, vol.~35, no.~1, pp. 93--109, jan 2015.

\bibitem{Satchidanandan2017}
B.~Satchidanandan and P.~R. Kumar, ``{Dynamic Watermarking: Active Defense of
  Networked Cyber–Physical Systems},'' \emph{Proc. IEEE}, vol. 105, no.~2,
  pp. 219--240, feb 2017.

\bibitem{Ding2018a}
\BIBentryALTinterwordspacing
D.~Ding, Q.~L. Han, Y.~Xiang, X.~Ge, and X.~M. Zhang, ``{A survey on security
  control and attack detection for industrial cyber-physical systems},''
  \emph{Neurocomputing}, vol. 275, pp. 1674--1683, 2018. [Online]. Available:
  \url{https://doi.org/10.1016/j.neucom.2017.10.009}
\BIBentrySTDinterwordspacing

\bibitem{Zhao2020}
\BIBentryALTinterwordspacing
Y.~Zhao and C.~Smidts, ``{A control-theoretic approach to detecting and
  distinguishing replay attacks from other anomalies in nuclear power
  plants},'' \emph{Prog. Nucl. Energy}, vol. 123, no. March, p. 103315, 2020.
  [Online]. Available: \url{https://doi.org/10.1016/j.pnucene.2020.103315}
\BIBentrySTDinterwordspacing

\bibitem{Hosseinzadeh2019}
M.~Hosseinzadeh, B.~Sinopoli, and E.~Garone, ``{Feasibility and Detection of
  Replay Attack in Networked Constrained Cyber-Physical Systems},'' \emph{2019
  57th Annu. Allert. Conf. Commun. Control. Comput. Allert. 2019}, pp.
  712--717, 2019.

\bibitem{Zhai2020}
L.~Zhai and K.~G. Vamvoudakis, ``{A data-based private learning framework for
  enhanced security against replay attacks in cyber-physical systems},''
  \emph{Int. J. Robust Nonlinear Control}, no. January, pp. 1--17, 2020.

\bibitem{Ferrari2017}
R.~M. Ferrari and A.~M. Teixeira, ``{Detection and Isolation of Replay Attacks
  through Sensor Watermarking},'' \emph{IFAC-PapersOnLine}, vol.~50, no.~1, pp.
  7363--7368, 2017.

\bibitem{Trapiello2019}
C.~Trapiello, D.~Rotondo, H.~Sanchez, and V.~Puig, ``{Detection of replay
  attacks in CPSs using observer-based signature compensation},'' \emph{2019
  6th Int. Conf. Control. Decis. Inf. Technol. CoDIT 2019}, pp. 1--6, 2019.

\bibitem{Sanchez2019}
H.~S. S{\'{a}}nchez, D.~Rotondo, T.~Escobet, V.~Puig, J.~Saludes, and
  J.~Quevedo, ``{Detection of replay attacks in cyber-physical systems using a
  frequency-based signature},'' \emph{J. Franklin Inst.}, vol. 356, no.~5, pp.
  2798--2824, 2019.

\bibitem{Ye2019}
\BIBentryALTinterwordspacing
D.~Ye, T.~Y. Zhang, and G.~Guo, ``{Stochastic coding detection scheme in
  cyber-physical systems against replay attack},'' \emph{Inf. Sci. (Ny).}, vol.
  481, no. 61773097, pp. 432--444, 2019. [Online]. Available:
  \url{https://doi.org/10.1016/j.ins.2018.12.091}
\BIBentrySTDinterwordspacing

\bibitem{Fang2020a}
C.~Fang, Y.~Qi, P.~Cheng, and W.~X. Zheng, ``{Optimal periodic watermarking
  schedule for replay attack detection in cyber–physical systems},''
  \emph{Automatica}, vol. 112, 2020.

\bibitem{Porter2020}
M.~Porter, P.~Hespanhol, A.~Aswani, M.~Johnson-Roberson, and R.~Vasudevan,
  ``{Detecting Generalized Replay Attacks via Time-Varying Dynamic
  Watermarking},'' \emph{IEEE Trans. Automat. Contr.}, vol.~66, no.~8, pp.
  1--1, 2020.

\bibitem{Rath2020}
S.~Rath, D.~Pal, P.~S. Sharma, and B.~K. Panigrahi, ``{A Cyber-Secure
  Distributed Control Architecture for Autonomous AC Microgrid},'' \emph{IEEE
  Syst. J.}, pp. 1--12, 2020.

\bibitem{Hoehn2016}
A.~Hoehn and P.~Zhang, ``{Detection of replay attacks in cyber-physical
  systems},'' \emph{Proc. Am. Control Conf.}, vol. 2016-July, pp. 290--295,
  2016.

\bibitem{Liu2020}
\BIBentryALTinterwordspacing
L.~Liu, L.~Ma, Y.~Wang, J.~Zhang, and Y.~Bo, ``{Distributed set-membership
  filtering for time-varying systems under constrained measurements and replay
  attacks},'' \emph{J. Franklin Inst.}, vol. 357, no.~8, pp. 4983--5003, 2020.
  [Online]. Available: \url{https://doi.org/10.1016/j.jfranklin.2020.01.029}
\BIBentrySTDinterwordspacing

\bibitem{Huang2020}
\BIBentryALTinterwordspacing
J.~Huang, L.~Zhao, and Q.~G. Wang, ``{Adaptive control of a class of strict
  feedback nonlinear systems under replay attacks},'' \emph{ISA Trans.}, no.
  xxxx, pp. 1--9, 2020. [Online]. Available:
  \url{https://doi.org/10.1016/j.isatra.2020.08.001}
\BIBentrySTDinterwordspacing

\bibitem{Gallo2018}
\BIBentryALTinterwordspacing
A.~J. Gallo, M.~S. Turan, F.~Boem, G.~Ferrari-Trecate, and T.~Parisini,
  ``{Distributed watermarking for secure control of microgrids under replay
  attacks},'' \emph{IFAC-PapersOnLine}, vol.~51, no.~23, pp. 182--187, 2018.
  [Online]. Available: \url{https://doi.org/10.1016/j.ifacol.2018.12.032}
\BIBentrySTDinterwordspacing

\bibitem{Franze2019}
G.~Franze, F.~Tedesco, and W.~Lucia, ``{Resilient Control for Cyber-Physical
  Systems Subject to Replay Attacks},'' \emph{IEEE Control Syst. Lett.},
  vol.~3, no.~4, pp. 984--989, 2019.

\bibitem{Chen2018b}
B.~Chen, D.~W. Ho, G.~Hu, and L.~Yu, ``{Secure Fusion Estimation for Bandwidth
  Constrained Cyber-Physical Systems under Replay Attacks},'' \emph{IEEE Trans.
  Cybern.}, vol.~48, no.~6, pp. 1862--1876, 2018.

\bibitem{Giraldo2019}
J.~Giraldo and A.~A. Cardenas, ``{A new metric to compare anomaly detection
  algorithms in cyber-physical systems},'' in \emph{Proc. 6th Annu. Symp. Hot
  Top. Sci. Secur.}, 2019, pp. 1--2.

\bibitem{Tartakovsky2014}
A.~Tartakovsky, I.~Nikiforov, and M.~Basseville, \emph{{Sequential analysis:
  Hypothesis testing and changepoint detection}}, 2014.

\bibitem{Girardin2018}
V.~Girardin, V.~Konev, and S.~Pergamenchtchikov, ``{Kullback-Leibler Approach
  to CUSUM Quickest Detection Rule for Markovian Time Series},'' \emph{Seq.
  Anal.}, vol.~37, no.~3, pp. 322--341, 2018.

\bibitem{watermarking_tac}
A.~names, ``{Quickest Detection of Deception Attacks in Networked Control
  Systems with Physical Watermarking, [other details to be added]},''
  \emph{arxiv org}, vol.~x, no.~x, pp. x--x, x.

\bibitem{Forsgren2002}
A.~Forsgren, P.~E. Gill, and M.~H. Wright, ``{Interior methods for nonlinear
  optimization},'' \emph{SIAM Rev.}, vol.~44, no.~4, pp. 525--597, 2002.

\bibitem{Boggs1995}
P.~T. Boggs and J.~W. Tolle, ``{Sequential Quadratic Programming},'' \emph{Acta
  Numer.}, vol.~4, no. 1995, pp. 1--51, 1995.

\bibitem{Johansson1998}
K.~H. Johansson and J.~L.~R. Nunes, ``{The Quadruple-Tank Process: A
  Multivariable Laboratory Process with an Adjustable Zero},'' \emph{Proc. Am.
  Control Conf.}, vol.~8, no.~3, pp. 456--465, may 2000.

\end{thebibliography}

\end{document}